\documentclass[a4paper,twoside,12pt]{amsart}

\usepackage[latin1]{inputenc}
\usepackage[english]{babel}

\usepackage{amsxtra}
\usepackage{amssymb}
\usepackage{amsmath,amsthm}
\usepackage[all]{xypic}
\usepackage{epsfig}
\xyoption{dvips}

\theoremstyle{definition}
\newtheorem{ntn}{Notation}[section]
\newtheorem{dfn}[ntn]{Definition}
\newtheorem{rem}[ntn]{Remark}
\newtheorem{exa}[ntn]{Example}
\theoremstyle{plain}
\newtheorem*{que}{Question}
\newtheorem*{thmintr}{Theorem}

\newtheorem{lem}[ntn]{Lemma}

\newtheorem{prp}[ntn]{Proposition}
\newtheorem{thm}[ntn]{Theorem}
\newtheorem{cor}[ntn]{Corollary}

\theoremstyle{remark}

\DeclareMathAlphabet{\mathds}{U}{dsrom}{m}{n}
\DeclareMathAlphabet{\mathsc}{U}{rsfs}{m}{n}

\DeclareMathOperator{\cha}{char}
\DeclareMathOperator{\id}{id}

\DeclareMathOperator{\Spec}{Spec}

\DeclareMathOperator{\Pic}{Pic}

\DeclareMathOperator{\NS}{NS}

\newcommand{\Q}{\mathbb{Q}}
\newcommand{\N}{\mathbb{N}}
\newcommand{\Z}{\mathbb{Z}}
\newcommand{\F}{\mathbb{F}}
\renewcommand{\P}{\mathbb{P}}
\newcommand{\R}{\mathbb{R}}

\newcommand{\C}{\mathbb{C}}

\renewcommand{\O}{\mathcal O}
\renewcommand{\L}{\mathcal L}

\newcommand{\M}{\mathcal M}
\newcommand{\E}{\mathcal E}

\newcommand{\D}{\mathcal D}
\newcommand{\X}{\mathcal X}

\renewcommand{\a}{\alpha}
\renewcommand{\l}{\lambda}

\newcommand{\e}{\epsilon}

\newcommand{\Fk}{\mathcal F}
\newcommand{\Av}{\mathcal A}

\newcommand{\lr}{\rightarrow}

\title{Non-Emptiness of the Height Strata of the Moduli Stack of Polarized K3 Surfaces}
\author{Jordan Rizov}
\address{Mathematisch Instituut\\ P.O. Box 80.010\\ 3508 TA Utrecht\\ The Netherlands}
\email{rizov@math.uu.nl}

\begin{document}
\maketitle
\begin{abstract}
In this paper we consider the following problem: For a given natural number
$d$ and a prime $p$ determine all Newton polygons of
polarized K3 surfaces of degree $2d$ over fields of characteristic $p$. This is
an analogue of the Manin problem for Newton polygons of abelian varieties. This question is equivalent to determining the non-empty height strata of the moduli stack $\M_{2d}\otimes \F_p$ of K3 surfaces with a polarization of degree $2d$ over $\F_p$. We prove here that if $d$ is large enough and prime to $p$, then the height strata of $\M_{2d}\otimes \F_p$ are non-empty.
\end{abstract}
\section*{Introduction}
In positive characteristic one can define interesting subvarieties of moduli
spaces of abelian varieties and of curves. Such loci can be given by considering the collection of these objects having fixed certain discrete invariants, such as for instance filtrations on ${\rm BT}_1$-groups or Newton polygons (see \cite{FO-Strat} and \cite{FO-EOStr}). A similar approach can be taken when studying moduli spaces of K3 surfaces. 

To every K3 surface over a perfect field $k$ of characteristic $p > 0$ one
associates a Newton polygon. By definition it is the Newton polygon of the
$F$-crystal $H^2_{\rm cris}(X/W(k))$. Denote by $\a$ the smallest slope of the Newton polygon of $X$. We define the height of $X$ to be infinite if $\a = 1$ and $1/(1-\a)$ otherwise. If
finite, the height of a K3 surface takes integral values from 1 to 10. Denote by $\M_{2d}$ the moduli stack of K3 spaces with a polarization of degree $2d$ and
suppose that $p$ does not divide $2d$. We look at the
subspaces $\M_{2d,\F_p}^{(h)}$ of $\M_{2d} \otimes \F_p$ of K3 surfaces with height at
least $h$. The collection of those 11 subspaces is called the height stratification of
$\M_{2d}\otimes \F_p$. One further stratifies $\M_{2d,\F_p}^{(11)}$ by the Artin
invariant (see for instance \cite{A}). In this way we obtain a filtration of the moduli space
$\M_{2d}\otimes \F_p$
\begin{displaymath}
\M_{2d}\otimes \F_p = \M_{2d, \F_p}^{(1)} \supset \M_{2d, \F_p}^{(2)} \supset
\dots \supset \M_{2d,\F_p}^{(11)} = \Sigma_1 \supset \dots \supset
\Sigma_{10}.
\end{displaymath}
The following question rises naturally.
\begin{que} 
Are all the subspaces in the height stratification of $\mathcal M_{2d}\otimes\F_p$ non-empty?
\end{que}
This can be reformulated in the following way: For a given natural number
$d$ and a prime $p$ determine all Newton polygons of
polarized K3 surfaces of degree $2d$ over fields of characteristic $p$. This is
an analogue of the Manin problem for  Newton polygons of abelian varieties
(\cite[Conj. 2, p. 76]{Man-ComGrSch}).

In this note we answer partially the question posed above by proving the following
result.
\begin{thmintr}
For every $d$, large enough and prime to $p > 2$, the subspaces in the height strata of $\M_{2d}\otimes \F_p$ are non-empty.
\end{thmintr}
The idea of the proof is to start with a polarized abelian surface $(A,\l)$ over
$k$ of certain degree and to use $\l$ to
construct an ample line bundle on the Kummer surface $X$ associated to $A$. In
this way we find $\bar k$-valued points of $\M_{2d} \otimes \F_p$. Making some appropriate
choices
of supersingular polarized abelian surfaces $(A,\l)$ we are able to show that
the height strata of $\M_{2d}\otimes \F_p$ are non-empty if $d$ is large enough. The
construction gives explicit bounds for $d$.

The organization of this note is the following. In Section \ref{HeightStratSection} we recall some definitions and give an overview of some results on the height stratification of $\M_{2d,\F_p}$. Section \ref{KummerSurfacesSection} is devoted to Kummer surfaces. Starting with an ample line bundle on an abelian surface we describe a way of constructing ample line bundles on its associated Kummer surfaces. This allows us to find points in $\M_{2d,\F_p}(\bar \F_p)$ which belong to certain height strata of $\M_{2d,\F_p}$. In Section \ref{KummerMapsSection} we take this idea one step further and construct Kummer morphisms from moduli stacks of polarized abelian surfaces to moduli stacks of polarized K3 surfaces. We use these morphisms to give an affirmative answer to the question posed above in case $d$ is large enough.
\newline
\newline
\newline
{\bf Notations}
\newline
\newline
We write $\M_{2d}$ for the Deligne-Mumford stack of K3 spaces with a polarization of degree $2d$. It is a smooth stack over $\Spec(\Z[1/2d])$. See \cite[\S 4.3, Thm. 4.7]{Riz-MK3} and \cite[Ch. 1, \S 1.4.3]{JR-Thesis}

If $A$ is a ring, $A \lr B$ a ring homomorphism then for any $A$-module ($A$-algebra etc.) $V$ we will denote by $V_B$ the $B$-module ($B$-algebra etc.) $V\otimes_A B$.

For an algebraic stack $\mathcal F$ over a scheme $S$ and a morphism of schemes $S' \lr S$ we will denote by $\mathcal F_{S'}$ the product $\mathcal F\times_S S'$ and consider it as an algebraic stack over $S'$.

We denote by $\Av_{g,d,n}$ the moduli stack of $g$-dimensional abelian varieties with a polarization of degree $d^2$ and a Jacobi level $n$-structure. It is a Deligne-Mumford stack which is smooth over $\Z[1/dn]$. We will write $\Av_{g,d}$ for $\Av_{g,d,1}$.
\newline
\newline
{\bf Acknowledgments}
\newline
\newline
This note contains the results of the second chapter of my Ph.D. thesis. I thank my advisors, Ben Moonen and Frans Oort for their help, their support and for everything I have learned from them. I thank the Dutch Organization for Research N.W.O. for the financial support with which my thesis was done.
\section{The Height Stratification of $\mathcal M_{2d,\F_p}$}\label{HeightStratSection}
Let $k$ be a perfect field of characteristic $p > 0$ and consider a K3 surface
$X$ over $k$. Consider the contravariant functor
\begin{displaymath}
 \Phi^2 : \underline{{\rm Art}} \lr {\rm Ab}
\end{displaymath}
from the category of local artinian schemes to abelian groups defined by
\begin{displaymath}
 \Phi^2(S) = \ker \bigl(H^2_{\rm et}(X \times S, \mathbb G_m) 
\lr H^2_{\rm et}(X, \mathbb G_m)\bigr).
\end{displaymath}
This functor is representable by a formal Lie group, denoted by $\widehat Br(X)$ and
called the \emph{formal Brauer group} of $X$.
\begin{prp} The formal group $\widehat Br(X)$ is $1$-dimensional and one has the following two
possibilities for it:
\begin{enumerate}
\item [(a)] The height of $\widehat Br(X)$ is infinite and then $\widehat Br(X) \cong {\widehat {\mathbb G}}_a$.
\item [(b)] The height is finite. Then $\widehat Br(X)$ is a $p$-divisible group. Moreover, its height satisfies $1 \leq h\bigl(\widehat Br(X)\bigr) \leq 10$.
\end{enumerate}
\end{prp}
For proofs we refer to \cite{A-M}. From now on we will call the
height of the formal Brauer group of $X$ simply the height of $X$
and denote it as $h(X)$. 

The \emph{Newton polygon} of $X$ is the Newton polygon of the $F$-crystal $H^2_{\rm cris}(X/W)$ where $W$ is the ring of Witt vectors $W(k)$ (see
\cite[\S 1.3 (c)]{Il-CrCoh}). It is the lower convex polygon starting at $(0,0)$ and ending at $(22,22)$. The height of a K3 surface $X$ can be read off from its Newton polygon. If $\a$ is the smallest slope of the Newton polygon of $X$, then $h(X) = 1/(1-\a)$ if $\a \ne 1$ and infinity otherwise. This follows from Corollary 3.3 in \cite[III]{A-M}.

The \emph{Hodge polygon in degree $m$} of a non-singular projective variety $X$ over $k$ is defined as the increasing convex polygon starting at $(0,0)$, having slope $i$ with multiplicity $h^{i,m-i} = \dim_k H^{m-i}(X,\Omega^i_X)$. For a K3 surface $X$ the Hodge polygon in degree $2$ will be called the Hodge polygon of $X$. The Newton polygon of a K3 surface lies on or above its Hodge polygon (\cite[Thm. 1.3.9]{Il-CrCoh}).

Recall that a K3 surface $X$ over $k$ is called \emph{ordinary} if any of the following equivalent conditions is satisfied:
\begin{enumerate}
\item [(i)] $h(X) = 1$.
\item [(ii)] The Newton and the Hodge polygon of $X$ coincide i.e., the Newton slopes of $X$ are 0 and 2 with multiplicity one, 1 with multiplicity 20.
\end{enumerate}
A K3 surface $X$ over $k$ is called \emph{supersingular} if any of the following equivalent conditions is satisfied:
\begin{enumerate}
\item [(i)] The height of $X$ is infinite.
\item [(ii)] The Newton polygon is a straight line i.e., all Newton slopes of $X$ are 1.
\end{enumerate}
The fact that the two possible ways of defining ordinary and supersingular K3 surfaces are equivalent follows from \cite[III, Cor. 3.3]{A-M}.

A K3 surface $X$ over $k$ is called \emph{supersingular in the sense of Shioda} if the rank of $\NS(X)$ is 22. One easily sees that if a K3 surface is
supersingular in the sense of Shioda, then it is supersingular. It is a
conjecture of M. Artin that, conversely, a supersingular K3 surface has N\'eron-Severi rank 22.
\begin{exa}
Let $p \equiv 3 \pmod 4$ be a prime number. Then the Fermat K3 surface
$$
 x^4 + y^4 + z^4 + w^4 = 0
$$
in $\P^3_{\F_p}$ is supersingular in the sense of Shioda (see \cite[Thm. 1]{Shi-CK3}).
\end{exa}
Let $d$ be an integer and assume further that $p$ does not divide $2d$. Consider the moduli stack $\mathcal M_{2d, \F_p} = \M_{2d} \otimes_\Z \F_p$ of K3 surfaces with
a polarization of degree $2d$ over a basis in characteristic $p$. Define the height stratification of $\M_{2d,\F_p}$ as follows: For $h \geq 1$ let $\M_{2d,\F_p}^{(h)}$ be the full subcategory of $\M_{2d,\F_p}$
\begin{displaymath}
 \M_{2d,\F_p}^{(h)}(S) = \bigl\{(X \lr S, \l) \in \M_{2d,\F_p}(S)\ |\ h(X_{\bar s})
\geq h\ {\rm for\ every\ geometric\ point}\ {\bar s} \in S \bigr\}.
\end{displaymath}
It is known that $\M_{2d,\F_p}^{(h)}$ is a closed substack of $\M_{2d,\F_p}$ of codimension at most $h - 1$. One defines a stratification of $\M_{2d,\F_p}^{(11)}$ by the Artin invariant (see \cite{A}). Let $X$ be a supersingular K3 surface and let $\Delta(\NS(X))$ be
the discriminant of the intersection pairing on $\NS(X)$
$$
 (\cdot ,\cdot ) \colon \NS(X) \times \NS(X) \lr \Z.
$$
One can show that ${\rm ord}_p(\Delta) = 2\sigma_0$ where $\sigma_0$ takes values $1, \dots, 10$. It is called
the \emph{Artin invariant} of $X$. Let $\Sigma_i$ be the full subcategory of $\M_{2d,\F_p}^{(11)}$ defined by
\begin{displaymath}
\begin{split}
 \Sigma_i(S) = \bigl\{(X \lr S,\l) \in \M_{2d,\F_p}(S)|& \ h(X_{\bar s}) = \infty\ {\rm and}\ \sigma_0(X_{\bar s})
\leq 11-i\\
&\ {\rm for\ every\ geometric\ point}\ {\bar s} \in S \bigr\}. \\
\end{split}
\end{displaymath}
In this way we obtain a filtration of the moduli space
\begin{equation}\label{StrFiltration}
\M_{2d,\F_p} = \M_{2d, \F_p}^{(1)} \supset \M_{2d, \F_p}^{(2)} \supset
\dots \supset \M_{2d,\F_p}^{(11)} = \Sigma_1 \supset \dots \supset
\Sigma_{10}.
\end{equation}
This is a chain of 20 closed substacks and the dimension drops with at least one at each step.
\begin{thm}\label{HStr}
For $h = 1,\dots ,10, 11$ the locus $\M_{2d,\F_p}^{(h)}$, if non-empty, is of codimension $h-1$ and for $h \ne 11$ is a local complete intersection.
\end{thm}
\begin{proof}
We refer to \cite[sect. 13,14 \& 15]{vdG-K}. The statement presented above is Theorem 15.1.
\end{proof}
\begin{rem}
B. Moonen and T. Wedhorn (\cite{MW-FZips}) have a theory of \emph{$F$-zips} which gives a scheme-theoretic and uniform definition of the filtration \eqref{StrFiltration}. For details we refer to Example 7.4 in {\it loc. cit.}.
\end{rem}
\section{Kummer Surfaces}\label{KummerSurfacesSection}
As we mentioned in the beginning of this chapter we will use Kummer surfaces to show that the height strata of $\M_{2d,\F_p}$ are non-empty for large enough $d$, prime to $p$. In Section \ref{KumSSection} we will recall some basic facts about Kummer surfaces which we will need in the sequel. In the next section, starting with a polarized abelian surface $(A,\l)$ we describe a way for constructing polarizations on its associated Kummer surface $X$. For this we make use of Seshadri constants.
\subsection{Kummer Surfaces}\label{KumSSection}
Recall that to an abelian surface $A$ over a field $k$ of characteristic
different from 2 we associated a K3 surface $X$, called the Kummer surface of $A$. We do this in the following way: Let $A[2]$ be the kernel of the multiplication by-2-map,
let $\pi \colon \tilde A \lr A$ be the blow-up of $A[2]$ and let $\tilde E$
be the exceptional divisor. The automorphism
$[-1]_A$ lifts to an involution $[-1]_{\tilde A}$ on $\tilde A$. Let $X$ be the quotient
variety of $\tilde A$ by the group of automorphisms $\{\id_{\tilde
A}, [-1]_{\tilde A}\}$ and denote by $\iota: \tilde A \lr X$ the
quotient morphism. It is a finite map of degree 2. The variety $X$ is the Kummer surface associated to $A$.

Assume further that all points in $A[2](\bar k)$ are $k$-rational. Then the exceptional divisor $\tilde E$ on $\tilde A$ consists of 16 irreducible curves ${E'}_j$, $j = 1,\dots, 16$, each corresponding to a point in $A[2](\bar k)$. We have that $({E'}_j,{E'}_l)_{\tilde A} = \delta_{j,l}$ where $\delta_{j,l}$ is the Kronecker $\delta$-function. Let us make the following notations:
\begin{displaymath}
\begin{split}   
{\E'}_j & := \O_{\tilde A}({E'}_j)\ \text{a\ line\ bundle\ on}\ \tilde A, \\
E_j & := \iota({E'_j})\ \text{a\ divisor\ on}\ X, \\
\E_j & := \O_X(E'_j)\ \text{the\ corresponding\ line\ bundle\ on}\ X.
\end{split}
\end{displaymath}
Then one has that
$$
 \iota^*\E_j \cong {\E_j'}^{\otimes 2}\ \ {\rm and}\ \ (\E_j,\E_l)_X = 2\delta_{j,l}.
$$
Moreover the line bundle
$\bigotimes_{j=1}^{16} \E_j$ is divisible by 2 in $\Pic(X_{\bar k})$. 

We turn next to some $p$-adic discrete invariants of Kummer surfaces. From now on we will assume that $k$ is a field of positive characteristic different from
$2$. Then $X$ is supersingular in the sense of Shioda if and only if $A$ is
supersingular. Indeed, according to \cite[\S 3, Prop. 3.1]{Shi-SSK3}, one has that $\NS(X_{\bar k})_\Q =
\NS(A_{\bar k})^{[-1]_A}_\Q \bigoplus_{j=1}^{16} \Q$. Hence  we have that $\text{rk}_\Z\NS(X) = 22$ if and only
if $\text{rk}_\Z\NS(A) =
6$ which is equivalent to $A$ being supersingular.
We will determine the Newton polygon of $X$ in term of the Newton polygon of $A$. To do that we shall need the following auxiliary result.
\begin{lem}\label{CohomKumS}
Let $A$ be an abelian surface and $X$ the associated Kummer surface over $k$. Then there
is a natural isomorphism
$$
 H_{\rm et}^2(X_{\bar k},\Q_l) \cong H_{\rm et}^2(A_{\bar k},\Q_l)
\oplus_1^{16} \Q_l(-1).
$$
\end{lem}
\begin{proof} This follows directly from the construction of Kummer surfaces. See for instance the proof of Lemma 2.2 in \cite{Ito-GRKS}.
\end{proof}
\begin{lem}\label{SrtIm}
Let $k$ be a finite field of characteristic different from 2. Then one has that
\begin{enumerate}
\item [(i)] If $A$ is ordinary, then the Newton polygon slopes of $X$ are
$\mu_1 =0;\ \mu_2 = \dots = \mu_5 = 1;\ \mu_6 = 2;\ \mu_j = 1$ for $j=7,\dots 22$.
In this case $X$ is ordinary i.e., its height is $1$.
\item [(ii)] If the $p$-rank of $A$ is 1, then the Newton polygon slopes of $X$ are
$\mu_1 = \mu_2 = 1/2;\ \mu_3 = \mu_4 = 1;\ \mu_5 = \mu_6 = 3/4;\  \mu_j = 1$ for $j=7,\dots 22$.
In this case $X$ has height is $2$.
\item [(iii)] If $A$ is supersingular, then $X$ is supersingular and all its Newton
polygon slopes are $1$. In this case $\Delta(\NS(A))$ and $\Delta(\NS(X))$ differ only by a power of $2$, hence the Artin invariant of $X$ is $1$ or $2$. It is $1$ if and only
if $A$ is a superspecial abelian surface.
\end{enumerate}
\end{lem}
\begin{proof}
As $k$ is a finite field and $X$ and $A$ are projective varieties one can compute the Newton polygons of $A$ and $X$ using \'etale cohomology instead of crystalline cohomology. We refer to \cite[1.3, Equality (1.3.5)]{Il-CrCoh} for an explanation and details. We will use the relation between the \'etale cohomology groups of $A$ and $X$ given in Lemma \ref{CohomKumS}.

If the Newton polygon slopes of $A$ are $\l_j$ for $j = 1,\dots, 4$, then those of $X$ satisfy
$\mu_1 = \l_1 + \l_2;\ \mu_2 = \l_1 + \l_3;\ \mu_3 = \l_1 + \l_4;\ \mu_4 = \l_2 + \l_3;\
\mu_5 = \l_2 + \l_4;\ \mu_6 = \l_3 + \l_4;\ \mu_i = 1$ for $i= 7, \dots, 22$.
The last statement follows from \cite[\S 3, Prop. 3.1]{Shi-SSK3}.
\end{proof}
\subsection{Ample Line Bundles on Kummer Surfaces}\label{AmpLBKumS}
Let $k$ be an algebraically closed field of characteristic different from 2 and consider an abelian surface $A$ over $k$. Denote by $X$ the associated Kummer surface. In this section we will show how to construct ample line bundles on $X$ starting with an ample bundle on $A$. This will allow us to give explicitly points in $\M_{2d,\F_p}$ for some $d$. 

Let $\L$ be an ample line bundle on $A$ with $\chi(\L) = d'$. Then by Riemann-Roch we have that $(\L,\L)_A = 2d'$. Let $n \in \N$ and fix 16 positive integers $n_j$. Consider the line bundle $\mathcal N$ on $\tilde A$ given by
$$
 \mathcal N = \pi^*\bigl(\L^n \otimes [-1]_A^*\L^n\bigr) \otimes
      \biggl(\bigotimes_{j=1}^{16} {\E_j'}^{-2n_j}\biggr).
$$
We will compute its self-intersection and show that $\mathcal N$ is the pull-back of a
line bundle on $X$.
\begin{lem}\label{extM}
 With the notations as above one has:
 \begin{enumerate}
 \item [(a)] $(\mathcal N,\mathcal N)_{\tilde A} = 8n^2d' - 4 \sum_{j=1}^{16}n_j^2$;
 \item[(b)]  There exists a line bundle $\M$ on $X$ such that $\iota^*\M \cong \mathcal N$. The line bundle $\M$ is ample iff $\mathcal N$ is ample. Moreover if $2d = (\M,\M)_X$, then we have that
 $$
  d = 2n^2d' - \sum_{j=1}^{16}n_j^2.
 $$
 \end{enumerate}
\end{lem}
\begin{proof}
(a) Combining \cite[Ch. V, \S 3, Prop. 3.2]{HAG} and the fact that $[-1]_A^*\L$ and
$\L$ are algebraically equivalent we get
\begin{displaymath}
\begin{split}
 (\mathcal N,\mathcal N)_{\tilde A} &= 2n^2(\L,\L)_A + 2n^2(\L,[-1]_A^*\L)_A - 4\sum_{j=1}{16}n_j^2 \\
 &= 8n^2d' - 4\sum_{j=1}{16}n_j^2.
\end{split}
\end{displaymath}
(b) Take a divisor $D \subset A \backslash A[2]$ such that $\O_A(D) = \L$. If $U := \tilde A \setminus \bigcup_{j = 1}^{16} E_j'$ and $V = X - \bigcup_{j = 1}^{16} E_j$, then the map $\iota \colon U \lr V$ is \'etale. Consider the divisor 
$$
 D_1 := \iota\bigl(n\pi^*(D) +n\pi^*([-1]_A^*D)\bigr)
$$ 
on $X$. As $D_1 \subset V$ we see that $\iota^*D_1 = n\pi^*(D) +n\pi^*([-1]_A^*D)$. Hence if we set 
$$
 \mathcal P := \O_X(D_1)
$$
then we have that $\iota^*\mathcal P \cong \pi^*\bigl(\L^n \otimes [-1]_A^*\L^n\bigr)$ on $\tilde A$. Using the fact that $\iota^*\E_j \cong {\E_j'}^{\otimes 2}$ one sees that the line bundle 
$$
 \M = \mathcal P \otimes \bigotimes_{j=1}^{16}\E_{j}^{\otimes - n_j}
$$
satisfies $\iota^*\M \cong \mathcal N$ on $\tilde A$. Since $\iota$ is a finite morphism $\M$ is ample on
$X$ if and only if $\mathcal N$ is ample on $\tilde A$ (\cite[Ch. III, Exercise 5.7 (d)]{HAG}).

For the self-intersection number computation one has
\begin{gather*}
 2n^2(\L,\L)_A + 2n^2(\L,[-1]_A^*\L) - 4 \sum_{j=1}^{16} n_j^2 = \\
 = (\mathcal N,\mathcal N)_{\tilde A}
 = (\iota^*\M,\iota^*\M)_{\tilde A} = \deg(\iota)(\M,\M)_X = 4d
\end{gather*}
which gives the formula from (b).
\end{proof}
\begin{rem}\label{graction}
Note that the line bundle $\L^n \otimes [-1]_A^*\L^n$ comes with a natural action of $[-1]_A$. Hence its pull-back $\pi^*\bigl(\L^n \otimes [-1]_A^*\L^n\bigr)$ comes equipped with an action of $[-1]_{\tilde A}$. Therefore one can apply \cite[Ch. III \S 10, Thm. 1(B)]{Mum-AV} to the morphism $\iota \colon \tilde A \lr X$ and conclude that $\mathcal P = \iota_*\bigl(\pi^*(\L^n\otimes [-1]_A^*\L^n)\bigr)^{[-1]_{\tilde A}}$ is the line bundle described in the proof of part (b).
\end{rem}
Lemma \ref{extM} suggests a way to construct ample line bundles on the Kummer surface
$X$. We will give sufficient conditions under which $\mathcal N$ is ample on $\tilde A$.
To do this we will make use of \emph{multiple Seshadri constants}. We will recall the definition below. For details we refer to \cite{Ba-S}.
\newline
\newline
{\bf Seshadri constants.} Let $\D$ be an ample line bundle on $A$ and let $x_1, \dots, x_{16}$ be the
points in $A[2](k)$ (recall that $k = \bar k$ and $\cha(k) \ne 2$). We make this change of notations here to avoid any possible confusion as later we will compute Seshadri constants for the ample line bundle $\D = \L \otimes [-1]^*_A \L^{-1}$. Let $\NS(\tilde A)_\R$ denote $\NS(\tilde A)\otimes_\Z \R$ and let $(\cdot , \cdot )_{\tilde A,\R}$ be the induced bilinear form. We will call an element $\mathcal R$ of $\NS(\tilde A)_\R$ \emph{numerically effective}, or shortly \emph{nef}, if for any irreducible curve $\Gamma$ in $\tilde A$ we have that $(\mathcal R,\O_{\tilde A}(\Gamma))_{\tilde A,\R} \geq 0$. Further, for an element $\mathcal R \in \NS(\tilde A)_\R$ and a real number $\e$ we will denote by $\mathcal R^\e$ the element $\e \cdot \mathcal R \in \NS(\tilde A)_\R$. 

One shows that
\begin{displaymath}
 \e_\D = \sup \biggl\{\e \in \R\ \bigg|\ \pi^*\D \otimes \bigotimes_{i=1}^{16}{\E_i'}^{-\e} \
{\rm is \ nef\ in\ } \NS(\tilde A)_\R \biggr\}
\end{displaymath}
exists. It is called the \emph{multiple Seshadri constant} on $A$ for
$x_1, \dots, x_{16}$. An equivalent definition of the Seshadri constant $\e$ can be given in the following way:
$$
 \e_\D = \inf \frac{(\D,\O_A(C))_A}
{\sum_{i=1}^{16}{\rm mult}_{x_i}C}
$$
where $\text{mult}_{x_i}C$ is the multiplicity of $C$ at $x_i$ and the infimum is taken over all irreducible curves $C$ in $A$ which pass through at least one $x_i$.
\begin{rem}\label{R_Sesh}
\begin{enumerate}
\item [(1)] If $0 < \delta < \e_\D$, then the line bundle
$\pi^*\D \otimes \bigotimes_{i=1}^{16}{\E_i'}^{-\delta}$ is nef. Moreover, one has the
strict inequality
$$
\bigl(\pi^*\D \otimes \bigotimes_{i=1}^{16}{\E_i'}^{-\delta}, \O_{\tilde A}(\Gamma)\bigr)_{\tilde A, \R} > 0
$$
for any irreducible curve $\Gamma$ on $\tilde A$.
\item [(2)] If $0 < n_i < \e_\D$, then $\pi^*\D \otimes \bigotimes_{i=1}^{16}{\E_i'}^{-n_i}$
is nef. Moreover, one has the strict inequality
$$
\bigl(\pi^*\D \otimes \bigotimes_{i=1}^{16}{\E_i'}^{-n_i}, \O_{\tilde A}(\Gamma)\bigr)_{\tilde A} > 0
$$
for any irreducible curve $\Gamma$ on $\tilde A$.
\end{enumerate}
These facts are clear from the second definition of $\e$.
\end{rem}
\noindent
{\bf Numerical estimates.} We will apply the general results on Seshadri constants to our particular situation. To avoid confusion let us make the following convention: If $A$ is an abelian surface, then by \emph{an elliptic curve $E$ in $A$} we shall mean an abelian subvariety $E$ of $A$ of dimension one.
\begin{prp}\label{Sesh}
Let $A$ be an abelian surface over an algebraically closed field $k$ of characteristic different from $2$. Let $\{x_1, \dots, x_{16}\}$ be the set of two-torsion points on $A$. Then for an ample line bundle $\D$ on $A$ we are in one of the following cases:
\begin{enumerate}
\item [(a)] The Seshadri constant satisfies the inequality
$$
 \e_\D \geq \frac{\sqrt {2(\D,\D)_A}} {16}.
$$
\item [(b)] The abelian surface $A$ contains a curve $E$ of genus $1$ such that
\begin{equation}\label{EllCSCost}
 \e_\D = \frac {(\D,\O_A(E))_A} {\#\{i | x_i \in E(k)\}}.
\end{equation}
\end{enumerate}
\end{prp}
\begin{proof}
See \cite[Prop. 8.3]{Ba-S}. Note that in this paper the assumption $k = \C$ is made.
However, the proof of the above proposition uses only the Hodge index theorem, the
Riemann-Roch theorem and some facts about blow-ups of curves. These results are valid over
any algebraically closed field.
\end{proof}
\begin{rem}
Note that in (b) we may assume that $E$ is an elliptic curve in $A$. Indeed, we have that $E$ is a translate of an elliptic curve $E' \subset A$ by a point $a \in A$. Since $\D$ is ample, the line bundles $t_a^*\D$ and $\D$ are numerically equivalent. Therefore we have that
$$
 (\D, \O_A(E)) = (t^*_a\D, t^*_a\O_A(E)) = (\D, t^*_a\O_A(E)) = (\D, \O_A(E')).
$$
We have further that $\#\{i | x_i \in E(k)\} \leq 4$. Indeed, all these points correspond to  points $p_i = x_i - a \in E'(k)$ for which $[2]p_i = [-2]a$ is a fixed point in $E'(k)$. As the isogeny $[2]$ is of degree $4$ (on $E'$) there are at most four such points. So we have that
\begin{displaymath}
\begin{split}
\e_\D &=  \frac {(\D,\O_A(E))_A} {\sum_{i=1}^{16}{\rm mult}_{x_i}E} = \frac {(\D,\O_A(E))_A} {\#\{i | x_i \in E(k)\}} = \\
 &=\frac {(\D,\O_A(E'))_A} {\#\{i | x_i \in E(k)\}} \geq \frac {(\D,\O_A(E'))_A}{4} = \frac {(\D,\O_A(E'))_A} {\sum_{i=1}^{16}{\rm mult}_{x_i}E'} \geq \e_\D. \\
\end{split}
\end{displaymath}
Therefore we have equalities and we conclude that $\#\{i | x_i \in E(k)\} = 4$. We also see that $a\in A[2](k)$.
\end{rem}

In what follows we will try to avoid case (b) of Proposition \ref{Sesh} as much as possible.
The reason is that one has little control over the intersection
$(\D,\O_A(E))_A$ in terms of the degree of $\D$. The bound in (a)
increases with $(\D,\D)_A$, but $A$ can contain curves of genus $1$ of any given intersection index $(\D,\O_A(E))_A$,
no matter how large $(\D,\D)_A$ is.

We will need the following auxiliary result which we shall apply to a line bundle $\L$ defining the polarization $\l$ on $A$ (cf. the beginning of this section).
\begin{lem}\label{Nak}
Let $\L$ be an ample line bundle on an abelian surface $A$ and let $E \subset A$ be an elliptic curve.
\begin{enumerate}
\item [(a)] Suppose that $(\L,\O_A(E))_A = 1$. Then there exists an elliptic curve
$E' \subset A$ such that $A \cong E \times E'$. Moreover, if $\pi_1 \colon E \times E' \lr E$
and $\pi_2 \colon E \times E' \lr E'$ are the two projections, then there exists a point
$P \in E$ and a line bundle $\mathcal G$ on $E'$ such that
$$
 \L \cong \pi_1^*\O_E(P) \otimes \pi_2^*{\mathcal G}.
$$
\item [(b)] Suppose that $(\L, \O_A(E))_A = m$ for some $m \in \N$. Then there exist an elliptic curve
$E' \subset A$ and an isogeny $f \colon E \times E' \lr A$ of degree at most $m$.
\end{enumerate}
\end{lem}
\begin{proof}
(a): The proof can be found in \cite[Lemma 2.6]{Nak}.

(b): Consider the homomorphism
$$
\xymatrix{
 \phi \colon A \ar[r]^{\varphi_\L} & A^t \ar[r] & E^t
}
$$
where $\varphi_\L$ is the map $a \mapsto t^*_a\L\otimes\L^{-1}$
and the second map is the dual of the inclusion $E \lr A$.
Let $E'$ be the reduced subscheme of the zero component of
$\text{ker}(\phi)$. Then $E'$ is an elliptic curve in $A$. Note that
$\L|_{E}$ is an invertible sheaf of degree at most $m$ hence $(E,E')_A \leq m$.
Define the homomorphism $E \times E' \lr A$ to be $(P,P') \mapsto P + P'$.
It is surjective and its kernel is a finite group scheme hence
it is an isogeny. Moreover its degree is exactly $(E,E')_A \leq m$.
\end{proof}
To get explicit conditions under which $\mathcal N$ is ample on $\tilde A$, one has to give some explicit estimates for $\e_\D$ for the ample line bundle $\D = \L^n \otimes [-1]^*_A\L^n$.
\begin{lem}\label{Ample}
With the notations of Lemma \ref{extM} one has that
\begin{enumerate}
\item [(a)] If $d',n,n_1, \dots n_{16}$ satisfy the following three inequalities
\begin{gather}
 2n^2d' - \sum_{i=1}^{16}n_i > 0\label{inequ1} \\
 n_i < \frac{n}{4}\label{inequ2} \\
 n_i < \frac{\sqrt{n^2d'}}{8}\label{inequ3},
\end{gather}
then the line bundle $\mathcal N$ is ample on $\tilde A$. 
\item [(b)] Assume further that $(A,\L)$ is not isomorphic to a polarized product of elliptic curves,
as in Lemma \ref{Nak} (a). Then for the ampleness of $\mathcal N$ on $\tilde A$ it is enough to require $n_i < n/2$ instead of \eqref{inequ2} along with the other two inequalities \eqref{inequ1} and~\eqref{inequ3}.
\end{enumerate}
\end{lem}
\begin{proof}
(a): Suppose that the inequalities \eqref{inequ1}, \eqref{inequ2} and \eqref{inequ3} are fulfilled. The first one simply says that $(\mathcal N,\mathcal N)_{\tilde A} > 0$. The second two are exactly the ones obtained from the explicit estimates for $\e_\D$.

Assume first that $\e_\D$ is computed by an elliptic curve $E$. Since $\L$ is ample on $A$ one has that $(\L,\O_A(E))_A \geq 1$. Hence by Proposition \ref{Sesh} we have that
$$
 \e_\D = \frac{\bigl(\L^n\otimes [-1]^*_A\L^n,\O_A(E)\bigr)_A}{4} = \frac{2n\bigl(\L, \O_A(E)\bigr)_A}{4}
\geq \frac{n}{2} \geq {2n_i}
$$
for every $i = 1, \dots, 16$.

If $\e_\D$ is not computed by by an elliptic curve, then case (a) of Proposition \ref{Sesh} and
the fact that $(\mathcal D, \mathcal D)_A = 8n^2d'$ give the estimate
$$
 \e_\D \geq \frac{\sqrt{n^2d'}} {4} \geq {2n_i}.
$$
for every $i = 1,\dots, 16$.

Thus if we impose these numerical conditions \eqref{inequ1}, \eqref{inequ2} and \eqref{inequ3} on $n$, $d'$ and $n_i$, then by Proposition \ref{Sesh} we have that $2n_i < \e_\D$. Hence by Remark \ref{R_Sesh} one has that
$(\mathcal N,\O_{\tilde A}(\Gamma))_{\tilde A} > 0$ for any irreducible curve $\Gamma$ on $\tilde A$.
Therefore by the Nakai-Moishezon criterion (\cite[Ch. V, \S 1, Thm. 1.10]{HAG}) the line bundle $\mathcal N$ is ample.

(b): Suppose that $(A, \L)$ is not isomorphic to a polarized product of elliptic curves, then $(\L, \O_A(E))_A \geq 2$. If $\e_\D$ is computed by an elliptic curve $E$ we have that 
\begin{displaymath}
 \e_\D \geq n \geq 2n_i\
\end{displaymath}
for all $i = 1,\dots, 16$. Otherwise, just like in (a) one has that 
\begin{displaymath} 
\e_\D \geq \frac{\sqrt{n^2d'}} {4} \geq {2n_i}
\end{displaymath}
for all $i$. Hence by the argument given in the proof of part (a) the line bundle $\mathcal N$ is ample.
\end{proof}

\section{Kummer Maps and Non-Emptiness of the Height Strata}\label{KummerMapsSection}
In the preceding section we gave a way to construct points in $\M_{2d}(\bar \F_p)$ starting with points in $\Av_{2,d'}(\bar \F_p)$ for some well-chosen integers $d$ and $d'$. Here we will show that this actually gives rise to morphisms between the stacks $\Av_{2,d',2,\F_p}$ and $\M_{2d,\F_p}$. We call these maps Kummer morphisms and we give their construction in detail in Section \ref{KummerMorphSect}. We will use them in Section \ref{NEHStrSect} to produce supersingular points in $\M_{2d,\F_p}(\bar \F_p)$ for $d$ large enough. In this way we will conclude that the height strata of $\M_{2d,\F_p}$ are non-empty for these $d$.
\subsection{The Kummer Morphisms}\label{KummerMorphSect}
We already saw that starting with an ample line bundle $\L$ on $A$ with $\chi(\L) = d'$ and
fixing integers $n, n_1, \dots, n_{16} > 0$ one produces a K3
surface $X$ and a line bundle $\M$ on it. This bundle is ample if further the
numerical conditions from Lemma \ref{Ample} are satisfied by $d',n, n_1, \dots, n_{16}$. 
It turns out that the resulting line bundle $\M$ depends only on the class of $\L$ in
$\NS(A)$. In other words, it depends only on the polarization $\l_\L$ defined by $\L$. Indeed, the construction
$$
 \L \mapsto \iota_*\bigl(\pi^*(\L\otimes [-1]_A^*\L)\bigr)^{[-1]_{\tilde A}} 
$$
gives a homomorphism of group schemes $h \colon \Pic_{A/k} \lr \Pic_{X/k}$ and since $\Pic^0_{X/k}$ is trivial we see that $h$ vanishes on $\Pic^0_{A/k}$. 

Suppose given numbers $n, d'$ and $n_1, \dots , n_{16}$ satisfying the
inequalities from Lemma \ref{Ample} (a). Then using the remark made above one shows that starting with a
polarized abelian surface $(A,\l)$ over an algebraically closed field $k$ one gets a polarized K3
surface $(X,\M)$. We will generalize this construction to a general base $S$. To
do so let us first try to find a more intrinsic way of constructing the line
bundle $\M$.

Let $\L$ be an ample line bundle on $A$ and let $\l = \varphi_\L$. The polarization defined by $\L \otimes [-1]_A^*\L$ is $2\l$. Let $\mathcal P$ be the
Poincar\'e sheaf on $A \times A^t$, where $A^t$ is the dual abelian surface. One has an isomorphism $[-1]_A^*\mathcal P \cong \mathcal P$. Then $\mathcal D = (\id_A \times \l)^* \mathcal P$ is a symmetric ample line bundle on
$A$ coming with an action of the group $\{\id_A,[-1]_A\}$. Moreover, the polarization $\varphi_{\mathcal D}$ is exactly $2\l$. 
 
The line bundles $\mathcal D$ and $\L \otimes [-1]_A^*\L$ are isomorphic. Indeed, consider the composition
\begin{displaymath}
\xymatrix{
A \ar[rr]^\Delta & & A\times A \ar[rr]^{\id_A \times \varphi_\L} & & A\times A^t.
}
\end{displaymath}  
where $\Delta \colon A \lr A\times A$ is the diagonal. By construction, $(\id_A \times \varphi_\L)^*\mathcal P$ is the Mumford bundle $\Lambda(\L)$ on $A \times A$, which pulls-back to $[2]^*\L \otimes \L^{-2}$ under $\Delta$. By Corollary 3 in \cite[Ch. II \S 6]{Mum-AV} we have that
$$
 [2]^*\L = \L^3\otimes [-1]_A^*\L.
$$
So we conclude that
$$
 \L\otimes [-1]^*_A\L = (\id_A\times \varphi_\L)^*\mathcal P.
$$
We will use the bundle $\mathcal D$ to generalize the construction given in Section \ref{AmpLBKumS} in relative settings.

We need to make another observation in order to be able to define Kummer morphisms. In the previous section we worked over an algebraically closed field $k$. Then we made use of points in $A[2](k)$ which give rise to some exceptional divisors on the blow-up surface $\tilde A$. We will carry out the same idea in the relative case. In order to be able to consider these exceptional divisors in general, for instance if the field $k$ is not algebraically closed, we will be working with abelian surface with level $2$-structure.

Let $(A \lr S, \l, \a) \in \Av_{2, d',2}(S)$ be a polarized abelian scheme over a base
scheme $S$ with a Jacobi level $2$-structure $\a$. Let $\mathcal P$ be the Poincar\'e bundle on
$A\times_S A^t$ where $A^t$ is the dual abelian scheme of $A$. Denote by
$\mathcal D$ the symmetric relatively ample line bundle $(\id_A \times
\l)^*\mathcal P$ on
$A$ (see \cite[Ch. 1, \S1, 1.6]{CF-AV}). Consider the blow-up $\tilde A$ of
$A$ at $A[2]$. Then the automorphism $[-1]_A$ extends to an
involution $[-1]_{\tilde A}$ on $\tilde A$. One forms the quotient $X$ of
$\tilde A$ by the finite automorphism group $\{\id_{\tilde A} , [-1]_{\tilde
A}\}$. Further we use the sheaf $\mathcal D$ to construct a polarization on $\X$.
We consider the sheaf
$$
 \mathcal N = \pi^*\mathcal D^n \otimes
      \biggl(\bigotimes_{j=1}^{16} {\E_j'}^{-2n_j}\biggr)
$$
on $\tilde A$ where $\mathcal E_j'$ are the 16 exceptional sheaves. One uses then \cite[Ch. III \S 10, Thm. 1(B)]{Mum-AV} to conclude
that $\mathcal N$ comes from a sheaf $\mathcal M$ on $X$ as in Proposition \ref{extM} and Remark \ref{graction}. Clearly this generalizes
the construction we considered over an algebraically closed field $k$. The sheaf $\mathcal M$ is then
fiberwise ample and hence $S$-ample by Lemma 1.10 in \cite{Riz-MK3} (Lemma 1.1.10 in \cite{JR-Thesis}). This $S$-ample line
bundle gives rise to a polarization of $X$. Isomorphisms of
polarized abelian schemes with a Jacobi level $2$-structure are sent to isomorphisms of polarized K3 schemes in a natural way. In this way
we get a morphism of stacks
\begin{displaymath}
 K_{n,n_1, \dots, n_{16}} \colon \Av_{2, d', 2} \lr \M_{2d, \Z[1/2]}
\end{displaymath}
sending an object $(A\lr S, \l, \a) \in \Av_{2, d', 2}$ to the object $(X,\M) \in \M_{2d, \Z[1/2]}$. We summarize this in the theorem below.
\begin{thm}\label{Kummer}
Let $n, d', n_1, \dots, n_{16} \in \N$ and assume that they satisfy the
numerical conditions \eqref{inequ1}, \eqref{inequ2} and \eqref{inequ3} of Lemma \ref{Ample}. Then there
exists a morphism of algebraic stacks
$$
 K_{n,n_1, \dots, n_{16}} \colon \Av_{2, d', 2} \lr \M_{2d,\Z[1/2]}
$$
where $d = 2n^2d' - \sum_{j=1}^{16}n_j^2$. The morphism sends a polarized abelian
surface, to its associated Kummer surface with an ample line bundle.
\end{thm}
\begin{dfn}
For any set of numbers $n, d', n_1, \dots, n_{16}$ satisfying the inequalities \eqref{inequ1}, \eqref{inequ2} and \eqref{inequ3} we will call the morphism $K_{n,n_1, \dots, n_{16}} \colon \Av_{2, d', 2} \lr \M_{2d,\Z[1/2]}$ constructed in Proposition \ref{Kummer} the \emph{Kummer morphism} (or \emph{Kummer map}) defined by $n, d', n_1, \dots, n_{16}$.
\end{dfn}
Recall that there are some weaker conditions (Lemma \ref{Ample} (b)) under which a polarized
abelian surface, which is not isomorphic to a polarized product of elliptic curves, gives a polarized Kummer surface. We will deal with this case now. One has a natural map
$$
 p \colon  \Av_{1,1,2} \times \Av_{1,d',2} \lr \Av_{2,d', 2}
$$
sending a pair of polarized elliptic curves to their polarized product as in Lemma \ref{Nak}.
Consider the open substack
\begin{displaymath}
 \mathcal U_{2,d',2} = \Av_{2,d',2} \setminus p(\Av_{1,1,2} \times \Av_{1,d',2})
\end{displaymath}
As we saw above one can construct a polarized Kummer surface out of any such abelian surface.
In the same lines one gets
\begin{prp}\label{kum}
Let $n, d', n_1, \dots, n_{16} \in \N$ satisfy the conditions
of Lemma \ref{Ample} (b). Then there exists a morphism of stacks
\begin{displaymath}
 K_{n,n_1,\dots, n_{16}} \colon \mathcal U_{2, d', 2} \lr \M_{2d,\Z[1/2]}
\end{displaymath}
as constructed in Theorem \ref{Kummer} where $d = 2n^2d' - \sum_{j=1}^{16}n_j^2$.
\end{prp}
\begin{proof}
$K_{n,n_1,\dots, n_{16}}$ maps a polarized abelian surface
to the polarized Kummer surface and this time one has to impose
the milder conditions of Lemma \ref{Ample} due to the fact that
the polarized products of elliptic curves are excluded.
\end{proof}
\begin{rem}\label{StratImage}
Let $(A,\l, \a)$ be an a polarized abelian surface over a finite field $k$ of
characteristic $p > 2$. Then using Lemma \ref{SrtIm} we see that the point
$K_{n,n_1,\dots, n_{16}}((A,\l, \a))$ in $\M_{2d, \F_p}(k)$ belongs to
\begin{enumerate}
\item [(i)] $\M_{2d,\F_p}^{(1)} \setminus \M_{2d,\F_p}^{(2)}$ if $A$ is ordinary;
\item [(ii)] $\M_{2d,\F_p}^{(2)} \setminus \M_{2d,\F_p}^{(3)}$ if the $p$ rank of
$A$ is 1;
\item [(iii)] $\Sigma_9 \setminus \Sigma_{10}$ if $A$ is supersingular but not
superspecial;
\item [(iv)] $\Sigma_{10}$ if $A$ is superspecial.
\end{enumerate}
\end{rem}
\subsection{Non-Emptiness of the Height Strata}\label{NEHStrSect}
Fix a prime number $p > 2$. We will prove here that the height strata of $\M_{2d,
\F_p}$ are non-empty for every large enough $d$ prime to $p$. The idea is to use
the Kummer maps and show that $\M_{2d, \F_p}$ contains a supersingular Kummer surface. 
Then by Theorem \ref{HStr} all strata are non-empty and so one has the claimed
dimensions.
\begin{thm} \label{NES}
For every large enough $d$ prime to $p$ the height strata of $\M_{2d,\F_p}$ are
non-empty.
\end{thm}
We will need the following result first.
\begin{lem}\label{res}
Every residue class modulo $2 \times 9^2$ can be represented by an integer of the form $\sum_{j=1}^{16} n_j^2$ with $1 \leq n_j \leq 4$ for all $j$.
\end{lem}
\begin{proof}
Explicit calculation.
\end{proof}
\begin{rem}
We believe that the statement of the preceding lemma remains
valid for all $n \geq 9$. In other words, all residues modulo $2n^2$ can be represented by an integer of the form $\sum_{j=1}^{16} n_j^2$ with $1 \leq n_j < \frac{n}{2}$. This is true for $n \in [9, 45]$.
\end{rem}
\noindent
{\it Proof of Theorem \ref{NES}.} First note that if for a given
$d$ there exist numbers $d'$ and $n_1, \dots, n_{16}$ giving a Kummer map
$$
 K_{d',n_1,\dots, n_{16}} \colon \mathcal U_{2, d', 2, \F_p} \lr \M_{2d,\F_p},
$$
as in Proposition \ref{kum}, then the height strata of $\M_{2d,\F_p}$ are non-empty. This follows
from Remark \ref{StratImage} as one can find a supersingular point in $\mathcal
U_{2, d', 2, \F_p}$.

Take $n=9$ and let $d' \geq 26$ so that the conditions of Lemma
\ref{Ample} (b) give $n_j \in [1,4]$. By Lemma \ref{res} we can
pick up 162 sets of numbers $(n_1, \dots , n_{16}), \ 1 \leq n_j
\leq 4$ which define Kummer maps as above and such that $F(n_1,
\dots , n_{16}) = \sum_{j=1}^{16} n_j^2$ gives all possible
resides modulo $2\times 9^2$. Hence the images of $\mathcal
U_{2, d', 2, \F_p}$ under those Kummer maps land in $\M_{2d, \F_p}$ where
$d = 2\times 9^2 d' - \sum_{j=1}^{16} n_j^2$. Using this set
of 162 sixteen-uples $(n_1,\dots,n_{16})$ and letting $d' \geq 26$ vary we can construct Kummer maps
for $\Fk_{2d, \F_p}$ for all $d \geq 2\times 9^2 \times 26 - 16 = 4196$.
Therefore by the remark we
started with the height strata of these moduli stacks are non-empty. This proves the assertion.
\qed
\newline
\newline
Using Kummer maps we saw that the height strata of $\M_{2d,\F_p}$ are non-empty if $d \geq 4196$.
On the other hand one has that the Fermat K3 surface
$$
 x^4 + y^4 + z^4 + w^4 = 0
$$
in $\P^3$, which is a Kummer surface by \cite[Thm. 1]{Shi-CK3}, is supersingular if
$p \equiv 3 \pmod 4$. Hence using explicit Kummer surfaces one can show the non-emptiness of
the height strata in lower polarization degrees. Using the same ideas as above we will ``cut some more moduli points'' of abelian surfaces in order
to improve the estimates in Lemma \ref{Ample}. In this way we will lower the bound for $d$.

First we will settle the case when $d$ is even. Let as before $A$ be an abelian surface and let $X$ be
the associated Kummer surface. The invertible sheaf $\sum_{j=1}^{16} \mathcal E_j$
is divisible by 2 in $\Pic(X)$. Hence $\sum_{j=1}^{16} \mathcal E_j'$ comes from a line bundle on $X$ modulo
2 torsion in $\Pic(\tilde A)$. Note that this torsion has to come from $\Pic^0(A)$. So it does not change
neither our constructions nor the intersection indexes we were dealing with.
Consider the following subset of $\Av_{2,d'}(\bar \F_p)$
\begin{displaymath}
\begin{split}
 U_{2,d'}^3  = \bigl\{(A,\l) \in \Av_{2,d'}(\bar \F_p) |&  {\rm there\ does\
 not\ exist\ an\ isogeny}\\
  & E\times E \lr A\ {\rm of\ degree}\ < 3 \bigr\}.
\end{split}
\end{displaymath}
For any $d'$ the supersingular locus of $\Av_{2,d', \F_p}$ remains
non-empty because we exclude only finitely many
points of it. Let $(A, \l) \in U_{2,d'}^3$ and let $\L$ be any ample line bundle on $A$ defining the polarization $\l$. Then by Lemma
\ref{Nak} we have that $(\L, \O_A(E))_A \geq 3$ for every elliptic
curve $E$ in $A$. Taking this into account and following the
proofs of Lemma \ref{Ample} and Theorem \ref{Kummer} one
constructs a Kummer map of sets
\begin{displaymath}
 K_{d'} \colon U_{2,d'}^3 \lr \M_{2d, \F_p}(\bar \F_p)
\end{displaymath}
where $n = n_1 = \dots = n_{16} = 1$, $d' \geq 32$ and $d = 2d' -
16$. Hence by Remark \ref{StratImage} we can conclude that
\begin{cor}
For all even $d \geq 48$ prime to $p$ the height strata of $\M_{2d, \F_p}$ are
non-empty.
\end{cor}
For the odd case we will construct Kummer maps with $n_1 = 1, n_2
= \dots = n_{16} = 2$. Define as before the set
\begin{displaymath}
\begin{split}
 U_{2,d'}^8  = \bigl\{(A,\l) \in \Av_{2,d'}(\bar \F_p) |& {\rm there\
 does\ not\ exist\ an\ isogeny}\\
  & E\times E \lr A\ {\rm of\ degree}\ < 8 \bigr\}
\end{split}
\end{displaymath}
Take a point $(A,\l) \in U_{2,d'}^8$ and let $\L$ be any ample line bundle on $A$ defining the polarization $\l$. Then according to Lemma \ref{Nak} (b) we have that $(\L,\O_A(E))_A \geq 9$ for all elliptic curves $E$ in $A$. Just as
above one constructs a Kummer map of sets
\begin{displaymath}
 K_{d'} \colon U_{2,d'}^8 \lr \M_{2d,\F_p}(\bar \F_p)
\end{displaymath}
where $n = 1, n_1 = 1, n_2 = \dots n_{16} = 2, d' \geq 512$ and $d
= 2d' - 15\times 4 - 1 = 2d' - 61$. Using these maps we obtain the following result.
\begin{cor}
For every odd $d \geq 963$ prime to $p$ the height strata of $\M_{2d, \F_p}$
are non-empty.
\end{cor}


\begin{thebibliography}{vdGK00}

\bibitem[AM77]{A-M}
M.~Artin and B.~Mazur.
\newblock {F}ormal {G}roups {A}rising from {A}lgebraic {V}arieties.
\newblock {\em Ann. scient. \'Ec. Norm. Sup.}, 10:87--132, 1977.

\bibitem[Art74]{A}
M.~Artin.
\newblock {S}upersingular {K}3 {S}urfaces.
\newblock {\em Ann. scient \'Ec. Norm. Sup.}, 7(4e):543--568, 1974.

\bibitem[Bau99]{Ba-S}
Th. Bauer.
\newblock {S}eshadri {C}onstants on {A}lgebraic {S}urfaces.
\newblock {\em Math. Ann.}, 313:547--583, 1999.

\bibitem[FC90]{CF-AV}
G.~Faltings and C-L. Chai.
\newblock {\em {D}egenerations of {A}belian {V}arieties}.
\newblock Number~22 in Ergebnisse der Mathematik und ihrer Grenzgebiete.
  Springer-Verlag, 1990.

\bibitem[Har77]{HAG}
R.~Hartshorne.
\newblock {\em {A}lgebraic {G}eometry}, volume~52 of {\em Graduate Texts in
  Mathematics}.
\newblock Springer-Verlag, New-York, 1977.

\bibitem[Ill95]{Il-CrCoh}
L.~Illusie.
\newblock {C}rystalline {C}ohomology.
\newblock In {\em in Motives, S. Kleiman, U. Jansen and J.-P. Serre eds.},
  volume~55 of {\em Proc. of Symp. in Pure Math.}, pages 43--70. AMS, 1995.

\bibitem[Ito01]{Ito-GRKS}
T.~Ito.
\newblock {G}ood {R}eduction of {K}ummer {S}urfaces.
\newblock {\em preprint}, pages 1--17, 2001.

\bibitem[Man63]{Man-ComGrSch}
Y.~I. Manin.
\newblock {T}heory of {C}ommutative {F}ormal {G}roups over {F}ields of {F}inite
  {C}haracteristic.
\newblock {\em Russian Math. Surv.}, 18(6):1--83, 1963.

\bibitem[Mum74]{Mum-AV}
D.~Mumford.
\newblock {\em {A}belian {V}arieties}.
\newblock Oxford University Press, New-York, 1974.

\bibitem[MW04]{MW-FZips}
B.~Moonen and T.~Wedhorn.
\newblock {D}iscrete {I}nvariants of {V}arieties in {P}ositive
  {C}haracteristic.
\newblock {\em International Math. Research Not.}, 72:3855--3903, 2004.

\bibitem[Nak96]{Nak}
M.~Nakamaye.
\newblock {S}eshadri {C}onstants on {A}belian {V}arieties.
\newblock {\em Amer. J. Math.}, 118:621--635, 1996.

\bibitem[Oor01a]{FO-Strat}
F.~Oort.
\newblock Newton polygon strata in the moduli space of abelian varieties.
\newblock In {\em in Moduli of abelian varieties (Texel Island, 1999)}, volume
  195 of {\em Progr. Math.}, pages 417--440. Birkh\"auser, Basel, 2001.

\bibitem[Oor01b]{FO-EOStr}
F.~Oort.
\newblock A stratification of the moduli space of abelian varieties.
\newblock In {\em in Moduli of abelian varieties (Texel Island, 1999)}, volume
  195 of {\em Progr. Math.}, pages 345--416. Birkh\"auser, Basel, 2001.

\bibitem[Riz05a]{JR-Thesis}
J.~Rizov.
\newblock {\em {M}oduli of {K}3 {S}urfaces and {A}belian {V}ariaties}.
\newblock PhD thesis, University of Utrecht, 2005.

\bibitem[Riz05b]{Riz-MK3}
J.~Rizov.
\newblock {M}oduli {S}tacks of {P}olarized {K}3 {S}urfaces in {M}ixed
  {C}haracteristic.
\newblock {\em preprint}, pages 1--37, 2005, math.AG/0506120.

\bibitem[Shi75]{Shi-CK3}
T.~Shioda.
\newblock {A}lgebraic {C}ycles on {C}ertain {K}3 {S}urfaces in {C}haracteristic
  $p$.
\newblock In {\em in Proc. Int. Conf. on Manifolds, Tokyo 1973}, pages
  357--364. Univ. of Tokyo Press, 1975.

\bibitem[Shi79]{Shi-SSK3}
T.~Shioda.
\newblock {S}upersingular {K}3 {S}urfaces.
\newblock In {\em in Algebraic Geometry, Proc. Summer Meeting, Univ. of
  Copenhagen}, Lecture Notes in Mathematics, pages 564--591. Springer-Verlag,
  1979.

\bibitem[vdGK00]{vdG-K}
G.~van~der Geer and T.~Katsura.
\newblock {O}n a {S}tratification of the {M}oduli of {K}3 {S}urfaces.
\newblock {\em J. Eur. Math. Soc. (JEMS)}, 2(3):259--290, 2000.

\end{thebibliography}
\end{document}